\documentclass{article}

\usepackage{amsfonts}
\usepackage{amssymb}
\usepackage{amsthm}
\usepackage{amsmath}
\usepackage{newlfont}
\usepackage{graphicx}
\usepackage[all,arc]{xy}
\usepackage{epsfig}

\usepackage{multicol}
\usepackage[T1]{fontenc}
\newtheorem {proposition}{Proposition}[section]
\newtheorem {theorem}{Theorem}[section]
\newtheorem {lemma}{Lemma}[section]

\newtheorem {exercise}{Exercise}
\newtheorem {corollary}{Corollary}[section]

\author{{\DJ{}or\dj{}e Barali\' {c}}\\ {\small Mathematical Institute SASA}\\[-2mm] {\small Belgrade, Serbia}
}

\title{Note on Totally Skew Embeddings of \\ Quasitoric Manifolds over Cube}
\date{}
\begin{document}
\maketitle

\begin{abstract} Skew embeddings are introduced by Ghomi and
Tabachnikov in \cite{Gho-Tab}. They are naturally related to
classical problems in topology, such as the generalized vector
field problem and the immersion problem for real projective
spaces. In recent paper \cite{nas}, {totaly skew} embeddings are
studied by using the topological obstruction theory. In the same
paper it is conjectured that for every $n$-dimensional, compact
smooth manifold $M^n$ $(n>1)$,
$$N(M^n)\leq 4n-2\alpha (n)+1,$$ where $N(M^n)$ is defined as the smallest dimension $N$ such that
there exists a {\em totally skew} embedding of a smooth manifold
$M^n$ in $\mathbb{R}^N$.

We prove that for every $n$, there is a quasitoric manifold
$Q^{2n}$ for which the orbit space of $T^n$ action is a cube $I^n$
and
$$N(Q^{2n})\geq 8n-4\alpha (n)+1.$$  Using the combinatorial properties of cohomology ring $H^* (Q^{2n}, \mathbb{Z}_2)$, we construct an
interesting general non-trivial example different from known
example of the product of complex projective spaces.
\end{abstract}

\renewcommand{\thefootnote}{}
\footnotetext{This research was supported by the Grant 174020 of
the Ministry for Education and Science of the Republic of Serbia.}

\section{Introduction}

The studying of skew embeddings was started by Ghomi and
Tabachnikov in \cite{Gho-Tab}. Recall, that two lines in an affine
space $\mathbb{R}^N$ are called {\em skew} if they are neither
parallel nor have a point in common or equivalently if their
affine span has dimension $3$.  More generally, affine subspaces
$U_1,\ldots, U_l$ of $\mathbb{R}^N$ are called {\em skew} if their
affine span has dimension $ \mathrm{dim}(U_1)+\cdots
+\mathrm{dim}(U_l)+l-1$, in particular a pair $U,V$ of affine
subspaces of $\mathbb{R}^N$ is skew if and only if each two lines
$p\subset U$ and $q\subset V$ are skew. An embedding $f :
M^n\rightarrow \mathbb{R}^{N}$ of a smooth manifold is called {\em
totally skew} if for each two distinct points $x,y\in M^n$ the
affine subspaces $df(T_xM)$ and $df(T_yM)$ of\, $\mathbb{R}^N$ are
skew. Define $N(M^n)$ as the smallest $N$ such that there exists a
totally skew embedding of $M^n$ into $\mathbb{R}^N$.

Ghomi and Tabachnikov established a surprising connection of
$N(M^n)$ with some classical invariants of geometry and topology.
For example they showed \cite[Theorem~1.4]{Gho-Tab} that the
problem of estimating $N(\mathbb{R}^n)$ is intimately related to
the generalized vector field problem and the immersion problem for
real projective spaces,  as exemplified by the inequality
$$
N(\mathbb{R}^n)\geq r(n) +n
$$
where $r(n)$ is the minimum $r$ such that the Whitney sum
$r\xi_{n-1}$ of $r$ copies of the canonical line bundle over
$\mathbb{R}P^{n-1}$ admits $n+1$ linearly independent continuous
cross-sections.

\medskip
Another example (\cite[Theorem~1.2]{Gho-Tab}) is the inequality
$$
N(S^n)\leq n + m(n) +1
$$
where $m(n)$ is an equally well-known function defined as the
smallest $m$ such that there exists a non-singular, symmetric
bilinear form $B : \mathbb{R}^{n+1}\times
\mathbb{R}^{n+1}\rightarrow \mathbb{R}^m$. As a consequence they
deduced the inequalities $N(S^n)\leq 3n+2$ and \\ $N(S^{2k+1})\leq
3(2k+1)+1$.

Up to our knowledge, the exact values of $N(M^n)$ are known only
for $N(\mathbb{R}^1) = 3$, $ N(S^1)=4$ and $N(\mathbb{R}^2)=6$.
General upper and lower bounds are given by Ghomi and Tabachnikov
inequality
\begin{equation}\label{eqn:lower-upper}
2n+2\leq N(M^n)\leq 4n+1.
\end{equation}

In the papers \cite{St-Tab}, \cite{GS-1} and \cite{GS-thesis} some
more general conditions with multiple regularity are studied.

\medskip

In the paper \cite{nas}, slightly different approach to the
invariant $N (M^n)$ is used. Using the topological obstruction
theory the lover bound is improved for various classes of
manifolds, such as projective spaces (both real and complex),
products of projective spaces, Grassmannians, etc. Stiefel-Whitney
classes are obstructions to totaly skew embeddings and it is shown
\cite[Proposition~1.]{nas} and \cite[Corollary~4.]{nas}
\begin{theorem} \label{skewteo}
If $k:=\max \{i\mid \overline{w}_i(M)\neq 0\}$ then
$$N(M)\geq 2n+2k+1.$$
\end{theorem}

In the same paper, some evidence in favor the conjecture
\cite[Conjecture~20.]{nas} $$N(M^n)\leq 4n-2\alpha (n)+1,$$ for
compact smooth manifold $M^n$ $(n>1)$, where $\alpha(n)$ is the
number of non-zero digits in the binary representation of $n$.
R.~Cohen \cite{Cohen} in 1985 resolved positively the famous {\em
Immersion Conjecture}, predicted that any compact smooth
$n$-manifold for $n>1$ can be immersed in
$\mathbb{R}^{2n-\alpha(n)}$.

\bigskip

Quasitoric manifolds are a class of manifolds with well understood
cohomology ring which is determined by Davis-Januszkiewicz formula
\cite[Theorem~4.14, Corollary 6.8]{DaJan}. Other topological
invariants could be computed from the formula, and we are
particulary interested in Stiefel-Whitney clases. In monography
\cite{BuPan} by Buchstaber and Panov there is a nice exposition on
quasitoric manifolds and its combinatorial and geometrical
properties.

Let $P$ be a simple polytope of dimension $n$ with $m$ facets and
$M$ a quasitoric manifold of dimension $2n$ over $P$. Let $v_j$
($\mathrm{deg} v_j=2, j=1, \dots, m$) be Poincar\'{e} dual of
codimension two invariant submanifold $M_j$ in $M^{2n}$, thus to
each facet $F_j$ we assign $v_j$ because the image of the
characteristic submanifold $M_j$ of the orbit projection
$M\rightarrow P$ is the facet $F_j$. The equivariant cohomology
ring $H^\ast_{T^n} (M; \mathbb{Z})=H^\ast (ET \times_{T^n} M)$ of
$M$ has the following ring structure:
$$H^\ast_{T^n} (M; \mathbb{Z})\simeq \mathbb{Z}[v_1, \dots,
v_m]/\mathcal{I},$$ where $\mathcal{I}$ is the Stanley-Reisner
(the face) ideal of polytope $P$ in the polynomial ring
$\mathbb{Z}[v_1, \dots, v_m]$.

Let $\pi : ET \times_{T} M \rightarrow B T$ be the natural
projection. The induced homomorphism $$\pi^\ast : H^\ast (B T)=
\mathbb{Z} [t_1, \dots, t_n]\rightarrow H^\ast (ET \times_{T^n}
M)=H^\ast_{T^n} (M; \mathbb{Z})$$ could be described by a $n\times
m$ matrix $\Lambda= (\lambda_1, \dots, \lambda_m)$, where
$\lambda_j \in \mathbb{Z}^n$ ( $j=1, \dots, m$ ) corresponds to
the generator of Lie algebra isotropy subgroup of characteristic
submanifold $M_j$. The matrix $\Lambda$ is called
\textit{characteristic matrix} of $M$. Put $\lambda_j=(\lambda_{1
j}, \dots, \lambda_{n j})^t \in \mathbb{Z}^n$. Then we have
$$\pi^\ast (t_i)=\sum_{j=1}^m \lambda_{i j} v_j$$ and let
$\mathcal{J}$ be the ideal in $\mathbb{Z}[v_1, \dots, v_m]$
generated by $\pi^\ast (t_i)$ for all $i=1, \dots, n$. The
ordinary cohomology of quasitoric manifolds has the following ring
structure: $$H^\ast (M)\simeq \mathbb{Z}[v_1, \dots,
v_m]/(\mathcal{I}+\mathcal{J}).$$ The Stiefel-Whitney class can be
described by the following  \textit{Davis-Januszkiewicz formula}:
$$\omega (M)=\imath^\ast \prod_{i=1}^m (1+v_i),$$ where $\imath$
is the inclusion $\imath: M \rightarrow ET \times_T M$ and
$\imath^\ast$ is the induced homomorphism.

\bigskip

In Section $2$ we describe one special quasitoric manifold $M_I$
over cube $I^n$ by matrix $\Lambda_{M_{I}}$. We describe its
cohomology ring and calculate the Stiefel-Whitney class.

Section $3$ is devoted to calculation of the Stiefel-Whitney class
of normal bundle using the smart manipulations of binomial
coefficients in cohomology ring (with $\mathbb{Z}_2$
coefficients). We calculate the obstruction to totally skew
embedding of manifold $M_I$ and get the main result of the paper.

\section{Quasitoric manifold over cube}

\subsection{Matrix $\Lambda_{M_{I}}$ and cube $I^n$}

A quasitoric manifold $M$ is described by two key objects: its
orbit polytope $P$ and the characteristic matrix $\Lambda$. Two
quasitoric manifolds over the same polytope, but with distinct
characteristic matrices are different, with non-isomorphic
cohomology rings. Although, polytope $P$ with its combinatorics
gives a lot of informations on manifold itself, the characteristic
matrix $\Lambda$ is necessary to understand important topological
invariants of the quasitoric manifold.

\begin{figure}[h!h!]
\centerline{\epsfig{figure=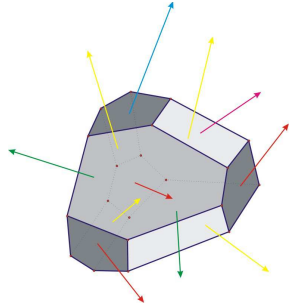, height=4 cm}}
\caption{Quasitoric manifold over polytope} \label{kvazpol}
\end{figure}

Let $\Lambda$ be the integer matrix $(n\times m)$ matrix whose
$i$-th column is formed by the coordinates of the facet vector
$\lambda_i$, $i=1, \dots, m$. Every vertex $v\in P^n$ is an
intersection of $n$ facets: $v=F_{i_1}\cap\cdots\cap F_{i_n}$. Let
$\Lambda_{(v)}:=\Lambda_{(i_1, \dots, i_n)}$ be the maximal minor
of $\Lambda$ formed by the columns $i_1$, $\dots$, $i_n$. Then
$$|\det \Lambda| =1.$$ In other words, to every facet $F_i$ there is
an assigned vector $\lambda_i=(\lambda_{1 i}, \dots, \lambda_{n
i})^t \in \mathbb{Z}^n$ in such  way that for every vertex
$v=F_{i_1}\cap\cdots\cap F_{i_n}$ vectors $\lambda_{i_1}, \dots,
\lambda_{i_1}$ are basis of lattice $\mathbb{Z}^n$ (see Figure
\ref{kvazpol})

\begin{figure}[h!h!]
\centerline{\epsfig{figure=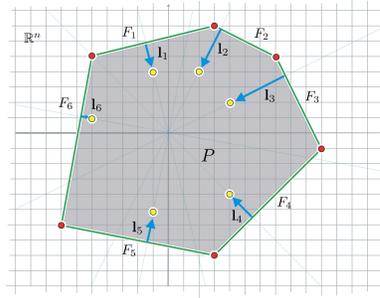, height=4 cm}}
\caption{Toric variety} \label{lepezaslika}
\end{figure}

In case when $P$ is a rational polytope and $\lambda_i\perp F_i$,
for every $i=1, \dots, m$ manifold $M$ is a toric variety (see
Figure \ref{lepezaslika}).

In monograph \cite[Construction~5.12]{BuPan} is described a
construction of quasitoric manifold from characteristic pair
$(P^n, l)$ that is from a combinatorial polytope $P$ and matrix
$\Lambda$.

\bigskip Let $I^n$ be a cube and $M_{I^n}$ a quasitoric manifold over $I^n$. Cube has $2 n$ facets $F_1$,
$\dots$, $F_n$, $F'_1$, $\dots$, $F'_n$ such that $F_i\cap
F'_i=\emptyset$ for every $i=1, \dots, n$. Let $v_1$, $\dots$,
$v_n$, $u_1$, $\dots$, $u_n$ be Poincar\'{e} duals to
characteristic submanifolds over the facets $F_1$, $\dots$, $F_n$,
$F'_1$, $\dots$, $F'_n$ respectively. Then Stanley-Reisner ideal
is generated by $$\mathcal{I}=\{v_1 u_1, v_2 u_2, \dots, v_n
u_n\}.$$

We study special quasitoric manifold $M_{I^n}$ over the cube, such
that vector $\lambda_j$ assigned to the facet $F_j$ (or the
generators of Lie algebra isotropy subgroup of characteristic
submanifold $M_j$) is $\lambda_j=(\underbrace{0, \dots, 0,}_{i-1}
1, \underbrace{0, \dots, 0}_{n-i})^t$ for every $j=1,\dots, n$ and
vector $\lambda_{j+n}$ assigned to the facet $F'_j$ is
$\lambda_{n+j}=(\underbrace{1, \dots, 1,}_{i} , \underbrace{0,
\dots, 0}_{n-i})^t$ for every $j=1,\dots, n$. Then we have:
$$\Lambda_{M_{I}}= \left[\begin{array}{cccccccc}
0 & 0 & \dots & 1 & 1 & 1 & \dots & 1\\
0 & 0 & \dots & 0 & 1 & 1 & \dots & 0\\
\vdots & \vdots & \ddots & \vdots & \vdots & \vdots & \ddots & \vdots\\
0 & 1 & \dots & 0 & 1 & 1 & \dots & 0\\
1 & 0 & \dots & 0 & 1 & 0 & \dots & 0\\
\end{array}\right].$$

Ideal $\mathcal{J}$ in $\mathbb{Z} [v_1, \dots, v_n, u_1, \dots,
u_n]$ is generated by linear forms $$\begin{array}{ccccccccc}
v_1 & + & u_1, &&&&&\\
v_2 & + & u_1 & + & u_2, &&&\\
\multicolumn{7}{c}{\dotfill},&\\
v_n & + & u_1 & + & u_2 & + & \cdots & + & u_n.
\end{array}$$

\subsection{Cohomology ring $H^\ast (M_{I})$ and Stiefel-Whitney
class $\omega(M_I)$}

Cohomology ring $H^\ast (M_{I})$ is determined using
Davis-Januszkiewicz theorem:

\begin{proposition}\label{kohomol} Cohomology ring $H^\ast (M_{I}; \mathbb{Z})$ is isomorphic
to $$H^\ast (M_{I}; \mathbb{Z})\simeq \mathbb{Z} [u_1, \dots,
u_n]/\mathcal{F}$$ where $\mathcal{F}$ is ideal in polynomial ring
$\mathbb{Z} [u_1, \dots, u_n]$ (such that $\mathrm{deg}
u_1=\dots=\mathrm{deg} u_n=2$) generated by quadratic forms
$$\mathcal{F}=\{u_1^2, u_2^2+u_1 u_2, \dots, u_n^2+u_1 u_n+u_2
u_n+\dots+u_{n-1} u_n\}.$$
\end{proposition}

It is easy to show the following relations in $H^\ast (M_{I};
\mathbb{Z})$:

\begin{proposition}\label{stgen} For every $i=1, \dots, n$ holds $$u_i^i\neq 0
\mbox{\,\, and\,\,} u_i^{i+1}=0.$$
\end{proposition}

\begin{proposition}\label{skrati} For every $i=2, \dots, n$ holds $$(1+u_i)(1+v_i)=1+u_1+\cdots+u_{i-1}$$
\end{proposition}

By Universal Coefficient Theorem we obtain that $$H^\ast (M_{I};
\mathbb{Z}_2)\simeq \mathbb{Z}_2 [u_1, \dots, u_n]/\mathcal{F}$$
where $\mathcal{F}$ is ideal in polynomial ring $\mathbb{Z}_2
[u_1, \dots, u_n]$ (such that $\mathrm{deg} u_1=\dots=\mathrm{deg}
u_n=2$) generated by quadratic forms
$$\mathcal{F}=\{u_1^2, u_2^2+u_1 u_2, \dots, u_n^2+u_1 u_n+u_2
u_n+\dots+u_{n-1} u_n\}.$$

Stiefel-Whitney class is the characteristic class in cohomology
with $\mathbb{Z}_2$ coefficients. By Davis-Januszkiewicz formula
Stiefel-Whitney class of $M_I$ is given by $$\omega
(M_I)=(1+u_1)\cdots (1+u_n)(1+v_1) \cdots (1+v_n),$$ but by using
Propositions \ref{kohomol} and \ref{skrati} it is easily reduced
to
$$\omega
(M_I)=(1+u_1)(1+u_1+u_2)\cdots (1+u_1+\cdots+u_{n-1}).$$

For the purposes of the main theorem, we are going to use another
form of cohomology ring $H^\ast (M_{I}; \mathbb{Z}_2)$. If we
choose another generators $t_1$, $\dots$, $t_n$ such that
$$\begin{array}{ccccccccc}
t_1 & = & u_1, &&&&&\\
t_2 & = & u_1 & + & u_2, &&&\\
\multicolumn{7}{c}{\dotfill},&\\
t_n & + & u_1 & + & u_2 & + & \cdots & + & u_n,
\end{array}$$ we get that $$H^\ast (M_{I};
\mathbb{Z}_2)\simeq \mathbb{Z}_2 [t_1, \dots, t_n]/\mathcal{G}$$
where $\mathcal{G}$ is ideal in polynomial ring $\mathbb{Z}_2
[t_1, \dots, t_n]$ (such that $\mathrm{deg} t_1=\dots=\mathrm{deg}
t_n=2$) generated by quadratic forms
$$\mathcal{G}=\{t_1^2, t_2^2+t_1 t_2, \dots, t_n^2+t_{n-1}
t_n\}.$$ Consequently, total Stiefel-Whitney class is given by
$$\omega (M_I)=(1+t_1) \cdots (1+t_{n-1}).$$

It is nod hard to check that the following is true in  $H^\ast
(M_{I}; \mathbb{Z}_2)$:

\begin{proposition}\label{stgen1} For every $i=1, \dots, n$ holds $$t_i^i=t_1 t_2\cdots t_i\neq 0
\mbox{\,\, and\,\,} t_i^{i+1}=0.$$
\end{proposition}

\section{Topological obstructions to totally skew embeddings of
manifold $M_I$}

\subsection{Stiefel-Whitney class $\overline{\omega} (M_I)$ of
normal bundle}

For the purposes of the main theorem, we are interested in
characteristic classes $\overline{\omega} (M_I)$ of the normal
bundle. Stiefel-Whitney classes  $ {\omega} (M_I)$ and
$\overline{\omega} (M_I)$ are related to each other by equality
$${\omega} (M_I)\cdot \overline{\omega} (M_I)=1.$$

In the previous section the Stiefel-Whitney classes  $ {\omega}
(M_I)$ is determined.  So, by Proposition \ref{stgen1}, it holds:

\begin{lemma}\label{swn} Total Stiefel-Whitney class $\overline{\omega} (M_I)$ of the normal
bundle is given by: $$\overline{\omega}
(M_I)=(1+t_1)(1+t_2+t_2^2)\cdots
(1+t_{n-1}+\cdots+t_{n-1}^{n-1}).$$
\end{lemma}

Since $\overline{\omega}_{2i} (M_I)=0$ when $i\geq n$, it is far
from obvious what is $\overline{\omega} (M_I)$ in cohomology ring
$H^\ast (M_{I}; \mathbb{Z}_2)$. For small $n$, we could calculate
$\overline{\omega} (M_I)$ by hand:

\begin{exercise}\begin{enumerate}
    \item $\overline{\omega} (M_{I^2}) = 1+t_1$,
    \item $\overline{\omega} (M_{I^3}) = 1+(t_1+t_2)$,
    \item $\overline{\omega} (M_{I^4}) = 1+(t_1+t_2+t_3)+t_1 t_3+t_1 t_2
    t_3$,
    \item $\overline{\omega} (M_{I^5}) = 1+(t_1+t_2+t_3+t_4)+(t_1 t_3+t_1 t_4+t_2 t_4)+(t_1 t_2 t_3+t_2 t_3
    t_4)$.
\end{enumerate}
\end{exercise}

By Lemma \ref{swn} for total Stiefel-Whitney classes of
$\overline{\omega} (M_{I^n})$ and $\overline{\omega}
(M_{I^{n+1}})$ the following recurrence relation holds (in $H^\ast
(M_{I^{n+1}}; \mathbb{Z}_2)$):
\begin{equation}\label{jed1}
\overline{\omega} (M_{I^{n+1}})=\overline{\omega} (M_{I^n})
(1+t_n+\cdots+t_n^n), \end{equation} or more explicitly
\begin{equation}\label{jed2}
\overline{\omega}_{2k} (M_{I^{n+1}})= \overline{\omega}_{2k}
(M_{I^n})+t_n\, \overline{\omega}_{2k-2} (M_{I^n})+\cdots+t_n^k
\,\mbox{for all}\, k=0, \dots, n-1
\end{equation} and \begin{equation}\label{jed3}
\overline{\omega}_{2n} (M_{I^{n+1}})= t_n
\,\overline{\omega}_{2n-2} (M_{I^n})+\cdots+t_n^n.\end{equation}

\medskip

Define numbers $\sigma^k_n$  for all positive integers $n$ and
$0\leq k \leq n-1$ as follows
\begin{equation}\label{sigmak}
\sigma^k_n=\left\{\begin{array}{rl} 1 & \mbox{if the total number
of the distinct monomials in\,} \overline{\omega}_{2k}
(M_{I^n}) \mbox{\, is odd}\\
0 & \mbox{\,elsewhere}\end{array}\right.
\end{equation}  So by \ref{jed2} and \ref{jed3}, it holds
\begin{equation}\label{sigmakrek}
\sigma^k_{n+1}=\sum_{i=0}^{k} \sigma^i_n
\end{equation} for every $k=1$, $\dots$, $n-1$ and
$$\sigma^{n}_{n+1}=\sigma^{n-1}_{n+1}.$$

Let us write the first $n$ rows of numbers $\sigma^k_n$ for $k=0$,
$\dots$, $n$:

$$\begin{array}{ccccccccc}
1 &  & &&&&&&\\
1 & 1 & & & &&&&\\
1 & 0 & 0& & &&&&\\
1 & 1 & 1 & 1 &&&&\\
1 & 0 & 1 & 0 & 0 &&&\\
1 & 1 & 0 & 0 & 0 & 0 &&\\
1 & 0 & 0 & 0 & 0 & 0 & 0 &\\
1 & 1 & 1 & 1 & 1 & 1 & 1 & 1\\
\multicolumn{9}{c}{\dotfill}
\end{array}$$

The previous sequence is closely related to the following sequence
of binomial coefficients ${n+k \choose k}$:
$$\begin{array}{ccccccccc}
\fboxsep 1 mm\fbox{\parbox{1.9 mm}{1}} &  & &&&&&&\\[1mm]
\fboxsep 1 mm\fbox{\parbox{1.9 mm }{1}} & \fboxsep 1 mm\fbox{\parbox{2.1 mm }{3}} & & & &&&&\\[2mm]
\fboxsep 1 mm\fbox{\parbox{1.9 mm}{1}} & 4 & 10& & &&&&\\[1mm]
\fboxsep 1 mm\fbox{\parbox{1.9 mm}{1}} & \fboxsep 1 mm\fbox{\parbox{2.1 mm }{5}} &  \fboxsep 1 mm\fbox{\parbox{3.5 mm }{$15$}} &  \fboxsep 1 mm\fbox{\parbox{3.5 mm }{$35$}} &&&&\\[2mm]
\fboxsep 1 mm\fbox{\parbox{1.9 mm}{1}} & 6 & \fboxsep 1 mm\fbox{\parbox{3.5 mm }{$21$}} & 56 & 70 &&&\\[2mm]
\fboxsep 1 mm\fbox{\parbox{1.9 mm}{1}} & \fboxsep 1 mm\fbox{\parbox{2.1 mm}{7}} & {8 \choose 2} & {9 \choose 3} & {10 \choose 4} & {11 \choose 5} &&\\[2mm]
\fboxsep 1 mm\fbox{\parbox{1.9 mm}{1}} & 8 & {9 \choose 2} & {10 \choose 3} & {11 \choose 4} & {12 \choose 5} & {13 \choose 6} &\\[2mm]
\fboxsep 1 mm\fbox{\parbox{1.9 mm}{1}} & \fboxsep 1
mm\fbox{\parbox{2 mm }{9}} &  \fboxsep 1 mm\fbox{\parbox{6 mm
}{${10 \choose 2}$}} & \fboxsep 1 mm\fbox{\parbox{6 mm }{${11
\choose 2}$}} & \fboxsep 1 mm\fbox{\parbox{6 mm }{${12 \choose
4}$}} & \fboxsep 1 mm\fbox{\parbox{6 mm }{${13 \choose 5}$}} &
\fboxsep 1 mm\fbox{\parbox{6 mm }{${14 \choose 6}$}} & \fboxsep 1
mm\fbox{\parbox{6 mm }{${15 \choose 7}$}}\\
\multicolumn{9}{c}{\dotfill}
\end{array}$$

Easy mathematical induction shows that:

\begin{lemma} $$\sigma^k_n\equiv {n+k \choose k} \pmod{2}.$$
\end{lemma}

By previous Lemma, in the case when $n=2^r$ we have
$$\sigma^{n-1}_n \equiv {2^r+ (2^r-1)\choose 2^r-1}\equiv {2^{r+1}-1\choose 2^r-1}\equiv 1
\pmod{2}.$$

Obviously, by the definition of $\sigma^k_n$ if $\sigma^k_n=1$
then $$\overline{\omega}_{2k} (M_{I^n})\neq 0.$$ Thus, we have:

\begin{theorem}\label{glavna1} If $n=2^r$ is a power of two then
$$\overline{\omega}_{2n-2} (M_{I^n})=t_1 t_2\cdots t_{n-1}\neq 0.$$
\end{theorem}

\begin{corollary}\label{posledica2nar} For $n=2^r$, quasitoric manifold
$M_{I^n}$ cannot be totally skew embedded in $\mathbb{R}^N$ when
$N$ is less than $$8 n - 3=4 \cdot\dim M_{I^n}-3.$$
\end{corollary}

\subsection{Topological obstructions when $n$ is not a power of $2$}

Theorem \ref{glavna1} is the sharpest possible result that one
could obtain using Stiefel-Whitney classes for quasitoric
manifolds. However, when $n$ is not a power of $2$ the previously
constructed quasitoric manifold $M_{I^n}$ in general does not
achieve the maximal possible value $k$ for which Stiefel-Whitney
class $\overline{\omega}_{2k} (M_{I^n})\neq 0$.

This problem could be overcome using the results from the previous
part.

\medskip

Let $n=2^{r_1}+2^{r^2}+\cdots+2^{r_t}$, $r_1>r_2>\cdots>r_t\geq 0$
be the binary representation of $n$ and let $m_i=2^{r_i}$ for
$i=1$, $\dots$, $t$. Divide the facets of a cube $I_n$ into $t$
groups $A_1$, $\cdots$, $A_t$ in such way that the opposite facets
belong to the same group and $|A_j|=2 m_j$ for $j=1$, $\dots$,
$t$. For every $j=1$, $\dots$, $t$,  denote the facets from $A_j$
with $F^{(j)}_i$ and $F'^{(j)}_i$ (the opposite facets), $i=1$,
$\dots$, $m_j$. We are going to construct new quasitoric manifold
$Q^{2n}$ over cube by defining a new characteristic matrix
$\Lambda$. Let $\lambda^{(j)}_i=(\underbrace{0, \dots,
0,}_{(\sum_{s}^{j-1} m_s)+(i-1)} 1, \underbrace{0, \dots,
0}_{n-(\sum_{s}^{j-1} m_s)-i})^t\in \mathbb{Z}^n$ and
$\lambda^{(j)}_{i+n}=(\underbrace{0, \dots, 0}_{(\sum_{s}^{j-1}
m_s)}\underbrace{1, \dots, 1,}_{i} \underbrace{0, \dots,
0}_{n-(\sum_{s}^{j-1} m_s)-i})^t\in \mathbb{Z}^n$ be the vectors
assigned to the facets $F^{(j)}_i$ and $F'^{(j)}_i$ respectively.
Let $v^{(j)}_i$ and $u^{(j)}_i$ be Poincar\'{e} duals to the
characteristic submanifolds over the facets $F^{(j)}_i$ and
$F'^{(j)}_i$ respectively, for all facets. Then the characteristic
matrix $\Lambda$ is:

$$\Lambda= \left[\begin{array}{c|c}
I_n & \begin{array}{ccc} \left[\begin{array}{c}0\end{array}\right]
&  \dots & \left[\begin{array}{cccc}
 1 & 1 & \dots & 1\\
 \vdots & \vdots & \ddots & \vdots\\
1 & 1 & \dots & 0\\
1 & 0 & \dots & 0
\end{array}\right]_{m_t\times m_t} \\
\left[\begin{array}{c}0\end{array}\right] &  \dots & \left[\begin{array}{c}0\end{array}\right] \\
\vdots &  \ddots & \vdots \\
\left[\begin{array}{c}0\end{array}\right] &
 \dots & \left[\begin{array}{c}0\end{array}\right] \\
\left[\begin{array}{cccc}
 1 & 1 & \dots & 1\\
 \vdots & \vdots & \ddots & \vdots\\
1 & 1 & \dots & 0\\
1 & 0 & \dots & 0
\end{array}\right]_{m_1\times m_1} &  \dots & \left[\begin{array}{c}0\end{array}\right]\\
\end{array}
\end{array}\right].$$

So by Davis-Januszkiewicz theorem we get:

\begin{theorem}\begin{itemize}
    \item Cohomology ring $H^\ast (Q; \mathbb{Z})$ is isomorphic
to $$H^\ast (Q; \mathbb{Z})\simeq \mathbb{Z} [u^{(1)}_1, \dots,
u^{(1)}_{m_1},\dots, u^{(t)}_{m_t}]/\mathcal{F}$$ where
$\mathcal{F}$ is ideal in polynomial ring $\mathbb{Z} [u^{(1)}_1,
\dots, u^{(1)}_{m_1},\dots, u^{(t)}_{m_t}]$ (such that
$\mathrm{deg} u^{(j)}_i=2$ for all $j$ and $i$) generated by
quadratic forms
$$\mathcal{F}=\{{u^{(j)}_1}^2, {u^{(j)}_2}^2+u^{(j)}_1 u^{(j)}_2, \dots, {u^{(j)}_{m_j}}^2+u^{(j)}_1 u^{(j)}_{m_j}+u^{(j)}_{2}
u^{(j)}_{m_j}+\dots+u^{(j)}_{m_j-1} u^{(j)}_{m_j} | j\in [t]\}.$$
    \item Cohomology ring $H^\ast (Q; \mathbb{Z}_2)$ is isomorphic
to $$H^\ast (Q; \mathbb{Z}_2)\simeq \mathbb{Z}_2 [u^{(1)}_1,
\dots, u^{(1)}_{m_1},\dots, u^{(t)}_{m_t}]/\mathcal{G}$$ where
$\mathcal{G}$ is ideal in polynomial ring $\mathbb{Z}_2
[u^{(1)}_1, \dots, u^{(1)}_{m_1},\dots, u^{(t)}_{m_t}]$ (such that
$\mathrm{deg} u^{(j)}_i=2$ for all $j$ and $i$) generated by
quadratic forms
$$\mathcal{G}=\{{u^{(j)}_1}^2, {u^{(j)}_2}^2+u^{(j)}_1 u^{(j)}_2, \dots, {u^{(j)}_{m_j}}^2+u^{(j)}_1 u^{(j)}_{m_j}+u^{(j)}_{2}
u^{j}_{m_j}+\dots+u^{(j)}_{m_j-1} u^{(j)}_{m_j} | j\in [t]\}.$$
    \item Total Stiefel-Whitney class $\omega(Q)$ is given by $$\omega
(Q)=\prod_{j=1}^t(1+u^{(j)}_{1})(1+u^{(j)}_{1}+u^{(j)}_{2})\cdots
(1+u^{(j)}_{1}+\cdots+u^{(j)}_{m_j-1}).$$
\end{itemize}
\end{theorem}

In the same fashion as in the previous section we choose new
generators $t_1^{(1)}, \dots, t_{m_1}^{(1)}, \dots, t_{m_t}^{(t)}$
such that $$H^\ast (Q; \mathbb{Z}_2)\simeq \mathbb{Z}_2
[t_1^{(1)}, \dots, t_{m_1}^{(1)}, \dots,
t_{m_t}^{(t)}]/\mathcal{G}$$ where $\mathcal{G}$ is ideal in
polynomial ring $\mathbb{Z}_2 [t_1^{(1)}, \dots, t_{m_1}^{(1)},
\dots, t_{m_t}^{(t)}]$ (such that $\mathrm{deg}
t_1=\dots=\mathrm{deg} t_n=2$) generated by quadratic forms
$$\mathcal{G}=\{{t^{(j)}_1}^2, {t^{(j)}_2}^2+t^{(j)}_1 t^{(j)}_2, \dots, {t^{(j)}_{m_j}}^2+t^{(j)}_{m_j-1}
t^{(j)}_{m_j}| j\in[t]\}.$$ Consequently, total Stiefel-Whitney
class is given by
$$\omega (Q)=\prod_{j=1}^t(1+t^{(j)}_{1}) \cdots (1+t^{(j)}_{m_j-1}).$$

Thus, the corresponding dual Stiefel-Whitney class is given by
$$\overline{\omega}
(Q)=\prod_{j=1}^t(1+t^{(j)}_{1})(1+t^{(j)}_{2}+{t^{(j)}_{2}}^2)
\cdots (1+t^{(j)}_{m_j-1}+\cdots+{t^{(j)}_{m_j-1}}^{m_j-1}).$$

But, according Theorem \ref{glavna1} we have:
$$\overline{\omega}
(Q)=\prod_{j=1}^t(1+(t^{(j)}_{1}+\cdots+{t^{(j)}_{m_j-1}})+\cdots+t^{(j)}_{1}t^{(j)}_{2}\cdots
t^{(j)}_{m_j-1}).$$

So we proved that the highest nontrivial dual Stiefel-Whitney
class is $$\overline{\omega}_{2n-2\alpha(n)} (Q)=t_1^{(1)}\cdots
t_{m_1-1}^{(1)}t_{1}^{(2)}\cdots t_{m_t-1}^{(t)},$$ where
$\alpha(n)$ is the number of non-zero digits in the binary
representation of $n$

As corollary we obtain:

\begin{theorem}[Main theorem] For every positive integer $n$ there is a
quasitoric manifold over cube $I^n$ such that $$N (Q)\geq 8
n-4\alpha(n)+1.$$
\end{theorem}

\emph{Remark.} Similar result cannot be obtained in the class of
toric varieties from a cube because Stiefel-Whitney class is
trivial in that case.

\vfill

{\small \DJ{}OR\DJ{}E BARALI\'{C}, Mathematical Institute SASA,
Kneza Mihaila 36, p.p.\ 367, 11001 Belgrade, Serbia

E-mail address: djbaralic@mi.sanu.ac.rs}

\end{document}